\documentclass[a4paper,11pt]{article}
\usepackage{amssymb,palatino}
\newtheorem{theorem}{Theorem}
\newtheorem{corollary}{Corollary}
\newtheorem{lemma}{Lemma}

\newcommand{\sty}{\displaystyle}
\newcommand{\QED}{\begin{flushright} $\Box$ \end{flushright}}
\newcommand{\Tr}{{\rm Tr} }                % traccia
\begin{document}

\title{Gauss Sums of the Cubic Character over $GF(2^m)$: \\ an elementary derivation}

\author{Davide Schipani\thanks{University of Zurich, Switzerland}, ~~
 Michele Elia \thanks{Politecnico di Torino, Italy}
}

\date{December 20, 2010}
\maketitle

\thispagestyle{empty}

\begin{abstract}
\noindent
An elementary approach is shown which derives the value of the Gauss sum of a cubic character over a finite field $\mathbb F_{2^s}$  without using Davenport-Hasse's theorem 
 (namely, if $s$ is odd the Gauss sum is $-1$, and if $s$ is even its value is $-(-2)^{s/2}$).
\end{abstract}


\vspace{2mm}

% ********************************************************************
\vspace{8mm}

\section{Introduction}

Let $\mathbb F_{2^{s}}$ be a Galois field over $\mathbb F_2$, and $\chi$ be the cubic character, namely $\chi$
 is a mapping from $\mathbb F_{2^{s}}^*$ into the complex numbers defined as
$$  \chi(\alpha^h \theta^j) = e^{\frac{2i\pi }{3}h} \dot{=} ~\omega^h ~~~~h=0,1,2~~,  $$ 
 where $\alpha$ is primitive and $\theta$ is a cube in $\mathbb F_{2^{s}}^*$, furthermore we set $\chi(0)=0$ by definition.

Let $\Tr_s(x)=\sum_{j=0}^{s-1} x^{2^j}$ be the trace function over $\mathbb F_{2^s}$, and $\Tr_{s/r}(x)=\sum_{j=0}^{s/r-1} x^{2^{rj}}$ be the relative trace function over $\mathbb F_{2^s}$
 relatively to $\mathbb F_{2^r}$, with $r|s$  \cite{lidl}. \\
%An Eisenstein sum of a character $\chi$ over $\mathbb F_{2^{s}}$ is defined as  \cite{berndt}
%$$ E_s(\chi)= \sum_{\stackrel{y \in  \mathbb F_{2^{s}}}{\Tr(y)=1}}  \chi(y)~~. $$ 
A Gauss sum of a character $\chi$ over $\mathbb F_{2^{s}}$ is defined as  \cite{berndt}
$$  G_s(\beta, \chi) = \sum_{y \in  \mathbb F_{2^s}}  \chi(y) e^{\pi i \Tr_s(\beta y)} = 
    \bar \chi(\beta) G_s(1, \chi) ~~~~\forall \beta \in \mathbb F_{2^{s}} ~~. $$

The values of the Gauss sums of a cubic character over $\mathbb F_{2^{s}}$ can be found by computing the Gauss sum over $GF(4)$ and applying Davenport-Hasse's theorem on the lifting of characters (\cite{berndt,Jungnickel,lidl}) for $s$ even (and by computing the Gauss sum over $GF(2)$ and then trivially lifting for $s$ odd). However it is possible to use a more elementary approach, and this is the topic of the present work.

If $s$ is odd then the cubic character is trivial because every element $\beta$ in $\mathbb F_{2^s}$ is a
cube as the following chain of equalities shows
$$ \beta \cdot 1 =  \beta \cdot (\beta^{2^s-1})^2= \beta \beta^{2^{s+1}-2} = \beta^{2^{s+1}-1} =(\beta^{\frac{2^{s+1}-1}{3}})^3 ~~,  $$   
since $\beta^{2^s-1}=1$, and $s+1$ is even, so that $2^{s+1}-1$ is divisible by $3$. In this case we have
$$ G_s(1, \chi) = \sum_{y \in  \mathbb F_{2^s}}  \chi(y) e^{\pi i \Tr_s( y)} = 
  \sum_{y \in  \mathbb F_{2^s}^*} e^{\pi i \Tr_s( y)} =-1  ~~, $$
since the number of elements with trace $1$ is equal to the number of elements with trace $0$ ($\Tr_s(x)\in\mathbb F_{2}$; moreover $\Tr_s(x)=1$ and $\Tr_s(x)=0$ are two equations of degree $2^{s-1}$), and
$e^{\pi i \cdot 0} =1$ while $e^{\pi i \cdot 1} =-1$. 

If $s=2m$ is even, the cubic character is nontrivial, and the computation of the Gauss sums requires some more effort; before we show how they can be computed with an elementary approach, we need some preparatory lemmas.

\section{Preliminary facts}

First of all we recall that, for any nontrivial character $\chi$ over $\mathbb F_{q}$, $\sum_{x \in \mathbb F_{q}} \chi(x) = 0 $.
This is used to prove a property of a sum of characters, already known to Kummer \cite{winter}, which can be
 formulated in the following form:
\begin{lemma}
  \label{lemma1}
Let $\chi$ be a nontrivial character and $\beta$ any element of $\mathbb F_{q}$; then
$$  \sum_{x \in \mathbb F_{q}} \chi(x) \bar \chi(x+\beta) = \left\{ \begin{array}{lcl} 
           q-1  &~~& \mbox{if} ~~\beta=0  \\ 
           -1      &~~& \mbox{if} ~~\beta\neq 0  
           \end{array} \right.   ~~. $$  
\end{lemma}

\noindent
{\sc Proof}.
If $\beta=0$, the summand is $\chi(x) \bar \chi(x) =1$, unless $x=0$ in which case it is $0$,
 then the conclusion is immediate. \\
When $\beta\neq 0$, we can exclude again the term with $x=0$, as $\chi(x)=0$, so that $x$ is invertible,
 and the summand can be written as 
$$\chi(x) \bar \chi(x+\beta)= \chi(x) \bar \chi(x) \bar \chi(1+\beta x^{-1}) = \bar \chi(1+\beta x^{-1}) ~~. $$
With the substitution $y=1+\beta x^{-1}$, the summation becomes
$$ \sum_{\stackrel{y \in \mathbb F_{q}}{y \neq 1}} \chi(y) = -1 + \sum_{y \in \mathbb F_{q}} \chi(y) =-1~~, $$
 as $\chi(y)=1$ for $y=1$. %\hfill $\Box$%
\QED

We are now interested in the sum $\sum_{x \in  \mathbb F_{q}}  \chi(x) \chi(x+1)$. Note that for the Gauss sums over $\mathbb F_{2^{s}}$ we have

\begin{equation}
  \label{gauss01}
G_{s}(1, \chi) = \sum_{\stackrel{y \in  \mathbb F_{2^{s}}}{\Tr_{s}(y)=0}}  \chi(y) -\sum_{\stackrel{y \in  \mathbb F_{2^{s}}}{\Tr_{s}(y)=1}}  \chi(y)~~.
\end{equation}
   
It follows that, if $\chi$ is a nontrivial character, then the Gauss sum over $\mathbb F_{2^{s}}$ satisfies the following: 
$$G_{s}(1,\chi) =2\sum_{\stackrel{y \in  \mathbb F_{2^{s}}}{\Tr_{s}(y)=0}}  \chi(y).$$ %=-2 E_{s}(\chi)
In fact 
 half of the field elements have trace $0$ and the other half $1$, so that  
$$\sum_{\stackrel{y \in  \mathbb F_{2^{s}}}{\Tr(y)=0}}  \chi(y) =-\sum_{\stackrel{y \in  \mathbb F_{2^{s}}}{\Tr(y)=1}}  \chi(y) $$
as the sum over all field elements is zero, since $\chi$ is nontrivial.%; the identity follows from (\ref{gauss01}). 
\begin{lemma}
   \label{lemma2}
If $\chi$ is a nontrivial character over $\mathbb F_{2^{s}}$, then
$$\sum_{x \in  \mathbb F_{2^{s}}}  \chi(x) \chi(x+1) = G_{s}(1,\chi)~~. $$   
\end{lemma}

\noindent
{\sc Proof}.
The sum  $\sum_{x \in  \mathbb F_{2^{s}}}  \chi(x) \chi(x+1)$ can be written as $\sum_{x \in  \mathbb F_{2^{s}}}  \chi(x(x+1))$, since the character is a
 multiplicative function, now the function $f(x)=x(x+1)$ is a mapping from $\mathbb F_{2^{s}}$ onto the subset of elements with trace $0$, as $\Tr_s(x)=\Tr_s(x^2)$ for any $s$, and each image comes exactly from two elements, $x$ and $x+1$. %, that have the same trace, since $\Tr_s(1)=0$ for $s$ even, which is our case. Therefore, half of the elements with trace $0$ are images of elements with trace $0$, and the remaining half are images of elements with trace $1$. 
It follows that 
\begin{equation}
   \label{gausssum}  
  \sum_{x \in  \mathbb F_{2^{s}}}  \chi(x) \chi(x+1) = 2 \sum_{\stackrel{y \in  \mathbb F_{2^{s}}}{\Tr_{s}(y)=0}}  \chi(y) = G_{s}(1,\chi) ~~.
\end{equation}        
\QED

\noindent
%%%Some properties of this function $G_m(1,\chi)$ are shown in the following theorems. \\

\begin{lemma}
   \label{lemma3}
Let $\chi$ be a nontrivial character of order $2^r+1$. Then the Gauss sum $G_{s}(1,\chi)$ is a real number. %, i.e. $G_{s}(1,\chi)=\pm 2^m$.
\end{lemma}

\noindent
{\sc Proof}. %Since $G_{2m}(1, \chi)$ is a sum of powers of $\omega$ the condition on its absolute value 
 %implies $G_{2m}(1, \chi)=(-1)^h \omega^j 2^{m}$, with $h\in \{0,1\}$ and $j\in \{0,1,2\}$. 
Using (\ref{gausssum})
 we have 
$$ \bar G_{s}(1,\chi)= \sum_{x \in  \mathbb F_{2^{s}}} \bar \chi(x) \bar \chi(x+1) = \sum_{x \in  \mathbb F_{2^{s}}} \chi(x^{2^r}) \chi(x^{2^r}+1) = \sum_{x \in  \mathbb F_{2^{s}}}  \chi(x)  \chi(x+1) = G_{s}(1,\chi) ~~,$$
as $\bar \chi(x)=\chi(x)^{2^r}=\chi(x^{2^r})$ and $x\to x^{2^r}$ is a field automorphism, so it just permutes the elements of the field. \hfill $\Box$

%Using (\ref{gausssum})
% we have 
%$$ \bar G_{2m}(1,\chi)= \sum_{x \in  \mathbb F_{2^{2m}}} \bar \chi(x) \bar \chi(x+1) = \sum_{x \in  \mathbb F_{2^{2m}}} \chi(x^2) \chi(x^2+1) = \sum_{x \in  \mathbb F_{2^{2m}}}  \chi(x)  \chi(x+1) = G_{2m}(1,\chi) ~~,$$
%because the conjugate of $\omega^j$ is $\omega^{2j}$, the character is multiplicative, and in $\mathbb F_{2^{2m}}$ the 
 %squaring operation is a field automorphism. \hfill $\Box$
%\QED 

\vspace{5mm}
\noindent 
\section{Main results}
The absolute value of $G_{s}(1,\chi)$ can be evaluated using elementary standard techniques going back to Gauss (see e.g. \cite{berndt}),
 while its argument requires a more subtle analysis.  Our main theorems in the following section derive in an elementary way the exact value of the Gauss sum for the cubic character $\chi$ over $\mathbb F_{2^{2m}}$ (the case of $s$ odd is trivial, as shown above). Before we proceed, we show in a standard way what is its absolute value.

Since $ G_{2m}(\beta, \chi)= \bar\chi(\beta) G_{2m}(1,\chi)$ , on one hand,  we have
\begin{equation}
  \label{sum1}
   \begin{array}{lcl}
  \sty   \sum_{\beta \in \mathbb F_{2^{2m}}} G_{2m}(\beta, \chi) \bar  G_{2m}(\beta, \chi) &=& \sty
     \sum_{\beta \in \mathbb F_{2^{2m}}}\bar \chi(\beta)  \chi(\beta) G_{2m}(1, \chi) \bar  G_{2m}(1, \chi) \\
      &=&  \sty \sum_{\beta \in \mathbb F_{2^{2m}}^*}  G_{2m}(1, \chi) \bar  G_{2m}(1, \chi) =
     (2^{2m}-1) G_{2m}(1, \chi) \bar  G_{2m}(1, \chi)~~. \end{array}
\end{equation}     
On the other hand, by the definition of Gauss sum, we have
$$   \begin{array}{l}
 \sty \sum_{\beta \in \mathbb F_{2^{2m}}} G_{2m}(\beta, \chi) \bar  G_{2m}(\beta, \chi) = \sum_{\beta \in \mathbb F_{2^{2m}}}\sum_{\alpha \in \mathbb F_{2^{2m}}}\sum_{\gamma \in \mathbb F_{2^{2m}}} \bar \chi(\alpha)
     e^{\pi i \Tr_{2m}(\beta \alpha)}
       \chi(\gamma) e^{-\pi i \Tr_{2m}(\gamma \beta)}  
\end{array} ~~,
$$     
and substituting $\alpha=\gamma+\theta$ in the last sum, we have
\begin{equation}
   \label{sum2}   
 \sum_{\beta \in \mathbb F_{2^{2m}}} G_{2m}(\beta, \chi) \bar  G_{2m}(\beta, \chi) = \sum_{\gamma \in \mathbb F_{2^{2m}}}\sum_{\theta \in \mathbb F_{2^{2m}}} \bar \chi(\gamma+\theta) 
        \chi(\gamma)\sum_{\beta \in \mathbb F_{2^{2m}}} e^{\pi i \Tr_{2m}(\beta \theta)} = 2^{2m}(2^{2m}-1)  ~~,
\end{equation}      
as the sum on $\beta$ is $2^{2m}$ if $\theta=0$ and is $0$ otherwise, since the values of the trace are equally distributed, as said above; consequently the sum over
 $\gamma$ is $2^{2m}-1$ times $2^{2m}$, as $\chi(0)=0$.
From the comparison of (\ref{sum1}) with (\ref{sum2}) we get $ G_{2m}(1, \chi) \bar  G_{2m}(1, \chi)=2^{2m}$, then $|G_{2m}(1, \chi)|=2^{m}$.

%The sign of $G_{2m}(1,\chi)$ is usually difficult to compute, since in our case 
Few initial values are $G_2(1,\chi)=2$, $G_4(1,\chi)=-4$, $G_6(1,\chi)=8$, $G_8(1,\chi)=-16$, and $G_{10}(1,\chi)=32$, so a reasonable guess is $G_{2m}(1,\chi)=-(-2)^{m}$.
This guess is correct as proved by the following theorems.

\begin{theorem}
   \label{signodd}
If $m$ is odd, the value of the Gauss sum $G_{2m}(1,\chi)$ is $2^m$.
\end{theorem}  

\noindent
{\sc Proof}. 
Let $\alpha$ a primitive cubic root of unity in $\mathbb F_{2^{2m}}$, then it is a root of $x^2+x+1$. In other words, a root 
 $\alpha$ of $x^2+x+1$, which does not belong to $\mathbb F_{2^{m}}$, as $m$ is odd, can be used to define a quadratic 
 extension of this field, i.e. $\mathbb F_{2^{2m}}$, and the elements of this extension can be represented 
  in the form $x+ \alpha y$, with $x,y \in  \mathbb F_{2^{m}}$. Furthermore, the two roots $\alpha$ and $1+\alpha$
  of $x^2+x+1$ are either fixed or exchanged by any Frobenius automorphism; in particular the automorphism
  $\sigma^m(x) = x^{2^m}$ necessarily exchange the two roots as it fixes precisely all the elements of 
  $\mathbb F_{2^{m}}$, while $\alpha$ does not belong to this field, so that  $\sigma^m(\alpha)\neq \alpha$. 
Now, a Gauss sum $G_{2m}(1,\chi)$ can be written as
\begin{equation}
   \label{signodd2}
 G_{2m}(1,\chi)= 2 \sum_{\stackrel{z \in  \mathbb F_{2^{2m}}}{\Tr_{2m}(z)=0}}  \chi(z) = 2 \sum_{\stackrel{x,y \in  \mathbb F_{2^{m}}}{\Tr_{2m}(x+\alpha y)=0}}  \chi(x+\alpha y) = 2 \sum_{\stackrel{x,y \in  \mathbb F_{2^{m}}}{\Tr_{m}(y)=0}}  \chi(x+\alpha y) ~~, 
\end{equation} 
 where we used the trace property 
$$\Tr_{2m}(x+\alpha y)=\Tr_{2m}(x)+\Tr_{2m}(\alpha y)=\Tr_{m}(x)+\Tr_{m}(x^{2^m})+\Tr_{2m}(\alpha y)=\Tr_{2m}(\alpha y),$$
 and the fact that 
$$ \begin{array}{lcl}
\Tr_{2m}(\alpha y) &=& \Tr_{m}(\alpha y)+\Tr_{m}(\alpha y)^{2^m}=\Tr_{m}(\alpha y)+\Tr_{m}((\alpha y)^{2^m}) \\
                   &=& \Tr_{m}(\alpha y)+\Tr_{m}(\alpha^{2^m} y) = \Tr_{m}(\alpha y)+\Tr_{m}((\alpha+1)y)=
                    \Tr_{m}(y)~~,  \end{array} $$
since $\alpha^{2^m}=\alpha+1$ as previously shown. 
The last summation in (\ref{signodd2}) can be split into three sums by separating the cases $x=0$ and $y=0$
$$  2 \sum_{\stackrel{x,y \in  \mathbb F_{2^{m}}}{\Tr_{m}(y)=0}}  \chi(x+\alpha y) = 
    2 \sum_{\stackrel{y \in  \mathbb F_{2^{m}}}{\Tr_{m}(y)=0}}  \chi(\alpha y) + 
    2 \sum_{x \in  \mathbb F_{2^{m}}}  \chi(x) + 
    2 \sum_{\stackrel{x,y \in  \mathbb F_{2^{m}}^*}{\Tr_{m}(y)=0}}  \chi(x+\alpha y) ~~.      $$ 
Considering the three sums separately, we have:
$$\sum_{x \in  \mathbb F_{2^{m}}}  \chi(x)=2^{m}-1 ~~,$$
 as $\chi(x)=1$ unless $x=0$ since $m$ is odd; 
$$\sum_{\stackrel{y \in  \mathbb F_{2^{m}}}{\Tr_{m}(y)=0}} \chi(\alpha y)= \chi(\alpha) (2^{m-1}-1) ~~,$$
 as the character is multiplicative, $\chi(y)=1$ unless $y=0$, and only the $0$-trace elements (which are $2^{m-1}-1$) should be counted;
$$  \sum_{\stackrel{x,y \in  \mathbb F_{2^{m}}^*}{\Tr_{m}(y)=0}}  \chi(x+\alpha y) =
     \sum_{\stackrel{x,y \in  \mathbb F_{2^{m}}^*}{\Tr_{m}(y)=0}} \chi(y) \chi(xy^{-1}+\alpha ) = 
     \sum_{\stackrel{ z,y \in  \mathbb F_{2^{m}}^*}{\Tr_{m}(y)=0}} \chi(z+\alpha ) =
    (2^{m-1}-1) \sum_{z \in  \mathbb F_{2^{m}}^*} \chi(z+\alpha )  ~~.    $$ 
as $y$ is invertible, $\chi(y)=1$ since $m$ is odd, $z$ has been substituted for $xy^{-1}$,
 and the sum we get in the end, being independent of $y$, is simply multiplied by the number of values assumed by
  $y$. Altogether we have
$$  G_{2m}(1,\chi)= 2^{m+1}-2+\chi(\alpha) (2^m-2) + (2^m-2)  \sum_{z \in  \mathbb F_{2^{m}}^*} \chi(z+\alpha ) =
      2^{m+1}-2 + (2^m-2)  \sum_{z \in  \mathbb F_{2^{m}}} \chi(z+\alpha )  ~~,    $$
      and, for later use, we define $A(\alpha) = \sum_{z \in  \mathbb F_{2^{m}}} \chi(z+\alpha )$. 
In order to evaluate $A(\alpha)$, we consider the sum of $A(\beta)$, 
 for every $\beta \in \mathbb F_{2^{2m}}$, and observe that $A(\beta)=2^m-1$ if $\beta \in \mathbb F_{2^{m}}$, while,
 if $\beta \not \in \mathbb F_{2^{m}}$ all sums assume the same value $A(\alpha)$, which is shown as follows: set $\beta=u+\alpha v$ with $v \neq 0$, then
$$ \sum_{z \in  \mathbb F_{2^{m}}} \chi(z+u+ \alpha v )  = \sum_{z \in  \mathbb F_{2^{m}}} \chi(v) \chi((z+u)v^{-1}+ \alpha )
    = \sum_{z' \in  \mathbb F_{2^{m}}}  \chi(z'+ \alpha )~~. $$
Therefore, the sum  $  \sum_{\beta \in  \mathbb F_{2^{2m}}} A(\beta) =\sum_{\beta}\sum_z \chi(z+\beta)=\sum_z\sum_{\beta} \chi(z+\beta)=0$ yields
$$  2^m(2^m-1)+(2^{2m}-2^m) A(\alpha) =0 $$
which implies $A(\alpha) =-1$, and finally
$$ G_{2m}(1,\chi)= 2^{m+1}-2 - (2^m-2) =2^m  ~~. $$ \hfill $\Box$

%\QED 
\paragraph{Remark 1.}
 
The above theorem can also be proved using a theorem by Stickelberger (\cite[Theorem 5.16]{lidl})

\begin{theorem}
   \label{signeven}
If $m$ is even, the Gauss sum $G_{2m}(1,\chi)$ is equal to $(-2)^{m/2} G_{m}(1,\chi)$ .
\end{theorem}  

\noindent
{\sc Proof}. 
The relative trace of the elements of $\mathbb F_{2^{2m}}$ over $\mathbb F_{2^{m}}$, which is 
$$   \Tr_{(2m/m)}(x)= x+x^{2^m} ~~,  $$
introduces the polynomial $x+x^{2^m}$ which defines a mapping from $\mathbb F_{2^{2m}}$ onto $\mathbb F_{2^{m}}$ with kernel the subfield $\mathbb F_{2^{m}}$ (\cite{lidl}). The equation $x^{2^m}+x=y$ has in fact exactly $2^m$ roots in $\mathbb F_{2^{2m}}$ for every $y \in \mathbb F_{2^{m}}$. \\
By definition we have
$$ G_{2m}(1,\chi)= 2\sum_{\stackrel{z \in  \mathbb F_{2^{2m}}}{\Tr_{2m}(z)=0}} \chi(z) =  
     2 \sum_{\stackrel{x,y \in  \mathbb F_{2^{m}}}{\Tr_{2m}(x+\alpha y)=0}}\chi(x+\alpha y)  ~~,  $$
where $\alpha$ is a root of an irreducible quadratic polynomial $x^2+x+b$ over $\mathbb F_{2^{m}}$, i.e. $\Tr_{m}(b)=1$ (\cite[Corollary 3.79]{lidl}) and $   \Tr_{(2m/m)}(\alpha)= 1$, which can be seen from the coefficient of $x$ of the polynomial. 
Now 
$$\Tr_{2m}(x+\alpha y)= \Tr_{2m}(x)+\Tr_{2m}(\alpha y)=\Tr_{2m}(\alpha y)= \Tr_{m}(\alpha y)+\Tr_{m}(\alpha^{2^m} y)~~,  $$
but $\alpha^{2^m}=1+\alpha$, so that $\Tr_{2m}(x+\alpha y)= \Tr_{m}(y)$, and we have
$$ G_{2m}(1,\chi)= 2\sum_{\stackrel{x,y \in  \mathbb F_{2^{m}}}{\Tr_{m}(y)=0}}\chi(x+\alpha y) =
 2\sum_{{x \in  \mathbb F_{2^{m}}}}\chi(x) + 2\sum_{\stackrel{y \in  \mathbb F_{2^{m}}^*}{\Tr_{m}(y)=0}}\chi(\alpha y)+
 2\sum_{\stackrel{x,y \in  \mathbb F_{2^{m}}^*}{\Tr_{m}(y)=0}}\chi(x+\alpha y)
 ~~,  $$
where the first summation has been split into the sum of three summations, by separating the cases $y=0$ and $x=0$.
We observe that, since the character over $\mathbb F_{2^{m}}$ is not trivial, the first sum is $0$ 
 and the second is $\chi(\alpha) G_{m}(1,\chi)$, while the third sum can be 
 written as follows
$$ 2\sum_{\stackrel{x,y \in  \mathbb F_{2^{m}}^*}{\Tr_{m}( y)=0}}\chi(x+\alpha y)= 
    2\sum_{\stackrel{x,y \in  \mathbb F_{2^{m}}^*}{\Tr_{m}( y)=0}} \chi(y) \chi(xy^{-1}+\alpha) =
    2\sum_{\stackrel{y \in  \mathbb F_{2^{m}}^*}{\Tr_{m}( y)=0}} \chi(y) 
    \sum_{z \in  \mathbb F_{2^{m}}^*}  \chi(z+\alpha)     
 ~~.  $$
 Putting all together, we obtain 
$$  G_{2m}(1,\chi)=  G_{m}(1,\chi) \sum_{z \in  \mathbb F_{2^{m}}}  \chi(z+\alpha) = G_{m}(1,\chi)  A_m(\alpha) ~~,  $$
which shows that  $| A_m(\alpha)|  = 2^{m/2}$ and that $A_m(\alpha)$ is real, as both $G_{2m}(1,\chi)$ and  $G_{m}(1,\chi)$ are real. Note that this holds for any $\alpha$ with $\Tr_{(2m/m)}(\alpha)=1$. \\

We will show now that 
 $A_m(\alpha)  = (-2)^{m/2}$. %To this aim, we firstly observe that $A_m(\alpha)$ is a real integer:
%$$ \bar A_m(\alpha) = \sum_{z \in  \mathbb F_{2^{m}}} \bar \chi(z+\alpha)=  \sum_{z \in  \mathbb F_{2^{m}}} \chi(z^2+\alpha^2) =
 %   \sum_{z \in  \mathbb F_{2^{m}}} \chi(z^2+\alpha+b) = \sum_{z' \in  \mathbb F_{2^{m}}} \chi(z'+\alpha) =
 %   A_m(\alpha) ~~. $$
%Secondly, if $\gamma \in \mathbb F_{2^{2m}}$ is any root of $q(x)=x^{2^m}+x+1$, the following argument shows that
% $A(\gamma)=A(\alpha)$: 
%observe preliminarily that $\alpha$ is a root of $q(x)$, and if $\gamma$ is another root of $q(x)$ then
 %$\gamma = u + \alpha$ for some $u \in \mathbb F_{2^{m}}$ , then we have 
%$$  A_m(\gamma) = \sum_{z \in  \mathbb F_{2^{m}}} \chi(z+\gamma) =
 %     \sum_{z \in  \mathbb F_{2^{m}}} \chi(z+u+\alpha) =
  %    \sum_{z' \in  \mathbb F_{2^{m}}} \chi(z'+\alpha) = A_m(\alpha) ~~. $$
Consider the sum of $A_m(\gamma)$ over all $\gamma$ with relative trace equal to $1$, which is, on one hand $2^m A_m(\alpha)$, as the polynomial $x^{2^m}+x=1$ has exactly $2^m$ roots in $\mathbb F_{2^{2m}}$ and on the other hand, explicitly we have
$$ \sum_{\stackrel{\gamma \in  \mathbb F_{2^{2m}}^*}{\Tr_{2m/m}(\gamma)=1}} A_m(\gamma) =
   \sum_{z \in  \mathbb F_{2^{m}}} \sum_{\stackrel{\gamma \in  \mathbb F_{2^{2m}}^*}{\Tr_{2m/m}(\gamma)=1}}  \chi(z+\gamma) =
   \sum_{z \in  \mathbb F_{2^{m}}} \sum_{\stackrel{\gamma' \in  \mathbb F_{2^{2m}}^*}{\Tr_{2m/m}(\gamma')=1}} \chi(\gamma')=
   2^m \sum_{\stackrel{\gamma' \in  \mathbb F_{2^{2m}}^*}{\Tr_{2m/m}(\gamma')=1}} \chi(\gamma')~~, $$
where the summation order has been exchanged, and $\Tr_{2m/m}(\gamma)=\Tr_{2m/m}(\gamma')$ as $\Tr_{2m/m}(z)=0$ for any $z\in\mathbb F_{2^{m}}$.
Comparing the two results, we have
$$  A_m(\alpha)= \sum_{\stackrel{\gamma' \in  \mathbb F_{2^{2m}}^*}{\Tr_{2m/m}(\gamma')=1}} \chi(\gamma')= 
     M_0+M_1 \omega+ M_2 \omega^2~~, $$
where $M_0$ is the number of $\gamma'$ with $\Tr_{2m/m}(\gamma')=1$ that are cubic residues, i.e. they have character
 $\chi(\gamma')$ equal to $1$, $M_1$ is the number of $\gamma'$ with $\Tr_{2m/m}(\gamma')=1$  that  have character 
 $\omega$, and $M_2$ is the number of $\gamma'$ with $\Tr_{2m/m}(\gamma')=1$ that have character $\omega^2$, then 
 $M_0+M_1+M_2=2^m$, and $M_1=M_2$ since $A_m(\alpha)$ is real. Therefore,we have $A_m(\alpha)= M_0-M_1$, and so we consider two equations for $M_0$ and $M_1$
$$   \left\{ \begin{array}{l}
       M_0 + 2 M_1 = 2^m  \\
       M_0-M_1 = \pm 2^{m/2}  \\
             \end{array}   \right.
$$
solving for $M_1$ we have $M_1= \frac{1}{3}(2^m \mp  2^{m/2}) $. Since $M_1$ must be an integer, we have
$$   \left\{ \begin{array}{lcl}
       M_0 - M_1 =  2^{m/2}    &~~& \mbox{if $m/2$ is even}\\
       M_0 - M_1 = -2^{m/2}&~~& \mbox{if $m/2$ is odd}.  \\
             \end{array}   \right.  ~~
$$

\QED
\begin{corollary}
If $m$ is even, the value of the Gauss sum $G_{2m}(1,\chi)$ is $-2^m$.
\end{corollary}
\noindent
{\sc Proof}.
It is a direct consequence of the two theorems above.
\QED

\section*{Acknowledgment}
The Research was supported in part by the Swiss National Science
Foundation under grant No. 126948

% ******************************************************************

\end{document}